\documentclass[11pt]{amsart}
\usepackage{latexsym}
\usepackage{amsmath}
\usepackage{amssymb}

\newcommand{\eproof}{\mbox{\ }\hfill $\Box$ \par \vskip 10pt}

\newtheorem{Theorem}{Theorem}[section]
\newtheorem{lemma}[Theorem]{Lemma}
\newtheorem{prop}[Theorem]{Proposition}

\newtheorem{corol}[Theorem]{Corollary}

\numberwithin{equation}{section}

\def\cal{\mathcal}

\begin{document}

\title[Semi-classical resolvent estimates]{Semi-classical resolvent estimates for $L^\infty$ potentials on Riemannian manifolds}

\author[G. Vodev]{Georgi Vodev}

\address {Universit\'e de Nantes, Laboratoire de Math\'ematiques Jean Leray, 2 rue de la Houssini\`ere, BP 92208, 44322 Nantes Cedex 03, France}
\email{Georgi.Vodev@univ-nantes.fr}

\date{}

\begin{abstract} We prove semi-classical resolvent estimates for the Schr\"odinger operator
with a real-valued $L^\infty$ potential on non-compact, connected Riemannian manifolds which may have a compact smooth boundary.
We show that the resolvent bound depends on the structure of the manifold at infinity. In particular, we show that for
compactly supported real-valued $L^\infty$ potentials and asymptoticaly Euclidean manifolds the resolvent bound is of the form 
$\exp(Ch^{-4/3}\log(h^{-1}))$, while for asymptoticaly hyperbolic manifolds it is of the form $\exp(Ch^{-4/3})$, where 
$C>0$ is some constant.
\end{abstract} 

\maketitle

\setcounter{section}{0}
\section{Introduction and statement of results}

The purpose of this paper is to extend the semi-classical resolvent estimates obtained recently in \cite{kn:KV}, \cite{kn:S2}, \cite{kn:V1} 
and \cite{kn:V2} for the Schr\"odinger operator in the Euclidean space $\mathbb{R}^n$ to a large class of non-compact, connected Riemannian manifolds $(M,g)$, $n={\rm dim}\, M\ge2$, with a smooth, compact
boundary $\partial M$ (which may be empty) and a smooth Riemannian metric $g$. We will consider manifolds of the form 
$M=X\cup Y$, where $X$ is a compact, connected Riemannian manifold with boundary $\partial X=\partial M\cup\partial Y$, while $Y$
is of the form $Y=[r_0,\infty)\times S$ with metric $g|_Y=dr^2+\sigma(r)$, where $(S,\sigma(r))$ is a compact $n-1$ dimensional
Riemannian manifold without boundary equipped with a family of Riemannian metrics $\sigma(r)$ depending smoothly on $r$ which can be written
in any local coordinates $\theta\in S$ in the form
$$\sigma(r)=\sum_{i,j}g_{ij}(r,\theta)d\theta_i d\theta_j,\quad g_{ij}\in C^\infty(Y).$$
Given any $r\ge r_0$, denote $Y_r=[r,\infty)\times S$. We can identify $\partial Y_r$ with the Riemannian manifold
$(S,\sigma(r))$. Then the negative Laplace-Beltrami operator on $\partial Y_r$ can be written in the form
$$\Delta_{\partial Y_r}=p^{-1}\sum_{i,j}\partial_{\theta_i}(pg^{ij}\partial_{\theta_j}),$$
where $(g^{ij})$ is the inverse matrix to $(g_{ij})$ and $p=(\det(g_{ij}))^{1/2}=(\det(g^{ij}))^{-1/2}$.
Let $\Delta_g$ denote the negative Laplace-Beltrami operator on $(M,g)$. Clearly, we can write the Laplace-Beltrami operator
$\Delta_Y:=\Delta_g|_Y$ in the form
$$\Delta_Y=p^{-1}\partial_r(p\partial_r)+\Delta_{\partial Y_r}=\partial_r^2+\frac{p'}{p}\partial_r+\Delta_{\partial Y_r}.$$
We have the identity
\begin{equation}\label{eq:1.1}
p^{1/2}\Delta_Yp^{-1/2}=\partial_r^2+\Lambda_\theta(r)-q(r,\theta)
\end{equation}
where 
$$\Lambda_\theta(r)=\sum_{i,j}\partial_{\theta_i}(g^{ij}(r,\theta)\partial_{\theta_j})$$
and $q$ is an effective potential given by the formula
$$q=(2p)^{-2}(\partial_rp)^2+(2p)^{-2}\sum_{i,j}g^{ij}\partial_{\theta_i}p\partial_{\theta_j}p-2^{-1}p\Delta_Y(p^{-1}).$$
We suppose that 
\begin{equation}\label{eq:1.2}
\sigma(r)\to f(r)^2\omega\quad\mbox{as}\quad r\to\infty
\end{equation}
where $\omega$ is a Riemannian metric on $S$ independent of $r$, which in the local coordinates $\theta\in S$ takes the form
$$\omega=\sum_{i,j}\omega_{ij}(\theta)d\theta_i d\theta_j,\quad \omega_{ij}\in C^\infty(S).$$
Here $f(r)$ is a function either of the form
\begin{equation}\label{eq:1.3}
f(r)=r^k,\quad k>0,
\end{equation}
or of the form 
\begin{equation}\label{eq:1.4}
f(r)=e^{r^\alpha},\quad 0<\alpha\le 1.
\end{equation}
The condition (\ref{eq:1.2}) implies
$$g^{ij}(r,\theta)\to f(r)^{-2}\omega^{ij}(\theta)\quad\mbox{as}\quad r\to\infty$$
where $(\omega^{ij})$ is the inverse matrix to $(\omega_{ij})$. In fact, we need stronger conditions on the functions $g^{ij}$,
namely the following ones:
\begin{equation}\label{eq:1.5}
\left|g^{ij}(r,\theta)-f(r)^{-2}\omega^{ij}(\theta)\right|\le Cf(r)^{-3},
\end{equation}
\begin{equation}\label{eq:1.6}
\left|\partial_r\left(g^{ij}(r,\theta)-f(r)^{-2}\omega^{ij}(\theta)\right)\right|\le Cf'(r)f(r)^{-4}
\end{equation}
with some constant $C>0$. Under the condition (\ref{eq:1.2}) we also have that the effective potential $q$ tends to the function
$$q_0(r)=\frac{(n-1)(n-3)f'(r)^2}{4f(r)^2}+\frac{(n-1)f''(r)}{2f(r)}.$$
More precisely, we suppose that for large $r$ the functions $q$ and $q_0$ satisfy
\begin{equation}\label{eq:1.7}
|q(r,\theta)-q_0(r)|\le Cr^{-1}f(r)^{-2},
\end{equation}
\begin{equation}\label{eq:1.8}
q_0(r)\le C,\quad \partial_rq_0(r)\le Cr^{-1}f(r)^{-2}
\end{equation}
with some constant $C>0$. In fact, an easy computation yields
$$q_0(r)=k(n-1)(kn-k-2)(2r)^{-2}$$
if $f$ is given by (\ref{eq:1.3}) and
$$q_0(r)=2^{-2}\alpha(n-1)(\alpha(n-1)+2(\alpha-1)r^{-\alpha})r^{2\alpha-2}$$
if $f$ is given by (\ref{eq:1.4}). Thus one can check that the condition (\ref{eq:1.8}) is always fulfilled 
if $f$ is given by (\ref{eq:1.4}), while in the other case it is fulfilled if $k\le 1$, $n\ge 2$, or
$k>1$, $n\ge 3$, or $k\ge 2$, $n=2$. In other words, (\ref{eq:1.8}) fails only in the case when $n=2$, $1<k<2$. 
 
Note that the above conditions are satisfied in the two most interesting cases which are
the asymptoticaly Euclidean manifolds (which corresponds to the choice $f(r)=r$) and 
the asymptoticaly hyperbolic manifolds (which corresponds to the choice $f(r)=e^r$).
In the first case we have $q_0=(n-1)(n-3)(2r)^{-2}$, while in the second case we have
$q_0=\left(\frac{n-1}{2}\right)^2$. 

Our goal is to study the resolvent of the Schr\"odinger operator
$$P(h)=-h^2\Delta_g+V(x)$$
where $0<h\ll 1$ is a semi-classical parameter and $V\in L^\infty(M)$ is a real-valued potential such that $V(r,\theta):=V|_Y$ satisfies the condition 
\begin{equation}\label{eq:1.9}
|V(r,\theta)|\le Cr^{-\delta}f(r)^{-2}
\end{equation}
with some constants $C>0$ and $\delta>1$. More precisely, we consider the self-adjoint realization of the operator $P(h)$
(which will be again denoted by $P(h)$) on the Hilbert space ${\cal H}=L^2(M,d{\rm Vol}_g)$. When the boundary $\partial M$ is not
empty we put Dirichlet boundary conditions. Given $s>1/2$ we let $\chi_s\in C^\infty(\overline{M})$, $\chi_s>0$, be a function
such that $\chi_s=1$ on $X$ and $\chi_s=r^{-s}$ on $Y_{r_0+1}$. 
We are going to bound from above the quantity
$$R_s^\pm(h,\varepsilon):=\log\left\|\chi_s(P(h)-E\pm i\varepsilon)^{-1}\chi_s
\right\|_{{\cal H}\to {\cal H}}$$
where $0<\varepsilon\le 1$ and $E>0$ is a fixed energy level independent of $h$. Set
$$m_0=
\left\{
\begin{array}{lll}
 \max\left\{\frac{2}{3k},\frac{1}{\delta-1}\right\}&\mbox{if $f$ is given by
 (\ref{eq:1.3})},\\
 \frac{1}{\delta-1}&\mbox{if $f$ is given by
 (\ref{eq:1.4})}.
\end{array}
\right.
$$
If $V$ is of compact support we set
$$m_0=
\left\{
\begin{array}{lll}
 \frac{2}{3k}&\mbox{if $f$ is given by
 (\ref{eq:1.3})},\\
 1&\mbox{if $f$ is given by
 (\ref{eq:1.4})}.
\end{array}
\right.
$$
Our main result is the following

\begin{Theorem} Let the potential $V$ satisfy (\ref{eq:1.9}). In the case when the function $f$ is given by
 (\ref{eq:1.4}) we suppose that $\delta>\frac{3\alpha}{4}+1$. 
Then there exist positive constants $C$ and $h_0$, independent
of $h$ and $\varepsilon$, such that for all $0<h\le h_0$ we have the bound
\begin{equation}\label{eq:1.10}
R_s^\pm(h,\varepsilon)\le 
 Ch^{-4/3-m_0(1-k)}\left(\log(h^{-1})\right)^{\frac{1-k}{\delta-1}}
\end{equation}
if $f$ is given by (\ref{eq:1.3}) with $k<1$. Moreover, if $V$ is of compact support we have the sharper bound
\begin{equation}\label{eq:1.11}
R_s^\pm(h,\varepsilon)\le Ch^{-\frac{2(k+1)}{3k}}.
\end{equation}
If $f$ is given by (\ref{eq:1.3}) with $k=1$ we have the bound 
\begin{equation}\label{eq:1.12}
R_s^\pm(h,\varepsilon)\le Ch^{-4/3}\log(h^{-1}).
\end{equation}
If $f$ is given by (\ref{eq:1.3}) with $k>1$ or by 
 (\ref{eq:1.4}) we have the bound 
\begin{equation}\label{eq:1.13}
R_s^\pm(h,\varepsilon)\le Ch^{-4/3}.
\end{equation}
\end{Theorem}

Recall that 
for asymptoticaly hyperbolic manifolds we have $f=e^r$, while for 
asymptoticaly Euclidean manifolds we have $f=r$. Thus we get the following

\begin{corol} Let $V\in L^\infty(M)$ be a compactly supported real-valued potential. Then, for asymptoticaly Euclidean manifolds of
dimension $n\ge 2$ we have the bound
(\ref{eq:1.12}), 
while for asymptoticaly hyperbolic manifolds of
dimension $n\ge 2$ we have the sharper bound (\ref{eq:1.13}).
\end{corol}

Note that for smooth potentials the following much sharper resolvent bound is known to hold (see \cite{kn:B2}, \cite{kn:D}, 
\cite{kn:S1})
\begin{equation}\label{eq:1.14}
R_s^\pm(h,\varepsilon)\le Ch^{-1}.
\end{equation}
 A high-frequency analog of (\ref{eq:1.14}) on Riemannian manifolds similar to the ones considered in the present paper 
 was also proved in 
\cite{kn:B1} and \cite{kn:CV}. In all these papers the regularity of the potential (and of the perturbation in general) 
is essential
in order to get (\ref{eq:1.14}). 
Without any regularity  the bound (\ref{eq:1.12}) has been
recently proved in \cite{kn:KV}, \cite{kn:S2} for real-valued compactly supported $L^\infty$ potentials 
when $M=\mathbb{R}^n$, $n\ge 2$, and in 
 \cite{kn:V1}, \cite{kn:V2} for real-valued short-range $L^\infty$ potentials 
when $M=\mathbb{R}^n$, $n\ge 3$. 
When $n=1$ it was shown in \cite{kn:DZ} that we have the better bound (\ref{eq:1.14}) instead of (\ref{eq:1.12}). 
When $n\ge 2$, however, the bound (\ref{eq:1.12}) seems hard to improve without extra conditions on the potential 
and it is not clear
if it is sharp or not. In contrast, it is well-known
that the bound (\ref{eq:1.14}) cannot be improved in general (e.g. see \cite{kn:DDZ}). 

To prove the above theorem we first prove a global uniform a priori estimate on an arbitrary compact, connected Riemannian manifold $X$. 
Roughly, we show that given an arbitrary open domain $U\subset X$, $U\neq\emptyset$, and any function $u\in H^2(X)$, we can control the Sobolev norm $\|u\|_{H^1(X)}$
by the norms $\|(P(h)-z)u\|_{L^2(X)}$ and $\|u\|_{H^1(U)}$, where $z\in\mathbb{C}$ (see Theorem 2.1 for the precise statement). When 
$\partial X$ is not empty we put Dirichlet boundary conditions on $u$. 
To do so we use the local Carleman estimates proved in \cite{kn:LR1} (see Proposition 2.2). We then propagate these local estimates in a way similar
to that one developed in \cite{kn:LR1}, making however some significant modifications due to the different nature of the problem we consider here.
 Note that local Carleman estimates with Neumann boundary conditions are proved in \cite{kn:LR2}, so most probably one can use the results in 
 \cite{kn:LR2} to conclude that Theorem 2.1 still holds in the case of Neumann boundary conditions. The proof, however, would be more
 technical and longer, and that is why we do not consider this case in the present paper.

In Section 4 we adapt the approach in \cite{kn:V1} to our situation in order to prove a global Carleman estimate on the end $Y$ of the
manifold $M$ (see Theorem 4.1). To do so, we construct in Section 3 global phase and weight functions on $Y$ in terms of
the function $f$, depending only on the
variable $r$, and we study their main properties. The most important one is the inequality (\ref{eq:3.9}) which is absolutely necessary 
for the Carleman estimate (\ref{eq:4.1}) to hold. Finally, in Section 5 we glue up the Carleman estimate on $Y$ with the a priori estimate on the
compact manifold $X$ comming from Theorem 2.1 to obtain the resolvent estimate. Note that a similar approach has already been used in
\cite{kn:KV} in the simpler case when $M=\mathbb{R}^n$, $n\ge 2$. In contrast, in \cite{kn:S2}, \cite{kn:V1} and \cite{kn:V2}
the global Carleman estimate is obtained on the whole space $\mathbb{R}^n$, which in turn poses some difficulties when 
$n=2$ due to the fact that in this case the effective potential (the function $q_0(r)$ above) is negative and the analysis
as $r\to 0$ gets quite complicated. That is why in \cite{kn:V1} and \cite{kn:V2} the condition $n\ge 3$ was imposed.
However, arguing as in the present paper we can avoid the problems related to the behaviour of the effective potential as $r\to 0$.
In fact, only the behaviour of the effective potential as $r\to \infty$ matters. 
Therefore the results in \cite{kn:V1} and \cite{kn:V2} hold for $n=2$, too.

\section{Carleman estimates on compact manifolds} 

Throughout this section $(X,g)$, $n={\rm dim}\,X\ge 2$, will be a compact, connected Riemannian manifold with a smooth boundary $\partial X$
which may be empty. Let $\Delta_g$ denote the negative Laplace-Beltrami operator on $(X,g)$ and introduce the operator
$$P(h)=-h^2\Delta_g+V(x)$$
where $0<h\le 1$ is a semi-classical parameter and $V\in L^\infty(X)$ is a complex-valued potential. 
Let $U\subset X$, $U\neq\emptyset$, be an arbitrary open domain, independent of $h$, such that $\partial U\cap\partial X=\emptyset$ and let $z\in\mathbb{C}$,
$|z|\le C_0$, $C_0>0$ being a constant independent of $h$. We will also denote by $H_h^1$ the Sobolev space equipped with the semi-classical norm. In this section we will prove the following

\begin{Theorem} There exists a positive constant $\gamma$ depending on $U$, $\sup|V|$ and $C_0$ but independent of $h$ such that 
for all $0<h\le 1$ we have the 
estimate
\begin{equation}\label{eq:2.1}
\|u\|_{H_h^1(X)}\le e^{\gamma h^{-4/3}}\|(P(h)-z)u\|_{L^2(X)}+e^{\gamma h^{-4/3}}\|u\|_{H_h^1(U)}
\end{equation}
for every $u\in H^2(X)$ such that $u|_{\partial X}=0$ if $\partial X$ is not empty.
\end{Theorem}

{\it Proof.} We will make use of the local Carleman estimates proved in \cite{kn:LR1}. Let $W\subset X$ be a small open domain and let
$x$ be local coordinates in $W$. If $\Gamma:=\overline W\cap\partial X$ is not empty we choose $x=(x_1,x')$, $x_1>0$ being the normal coordinate in $W$ and $x'$ the tangential ones. Thus in these coordinates $\Gamma$ is given by $\{x_1=0\}$. Let $p(x,\xi)\in C^\infty(T^*W)$ be the principal symbol of the operator $-\Delta_g$ and let $0<\hbar\ll 1$ be a new semi-classical parameter. Let $\varphi\in C^\infty(\overline W)$
be a real-valued function independent of $\hbar$. Then the principal symbol, $p_\varphi$, of the operator $-\hbar^2e^{\varphi/\hbar}\Delta_ge^{-\varphi/\hbar}$
is given by the formula
$$p_\varphi(x,\xi)=p(x,\xi+i\nabla\varphi(x)).$$
We suppose that $\varphi$ satisfies the H\"ormander condition
\begin{equation}\label{eq:2.2}
\forall (x,\xi)\in T^*W, p_\varphi(x,\xi)=0 \Longrightarrow \left\{{\rm Re}\,p_\varphi, {\rm Im}\,p_\varphi\right\}(x,\xi)>0.
\end{equation}
It is easy to check that (\ref{eq:2.2}) is fulfilled if we take $\varphi=e^{\lambda\psi}$, where $\psi\in C^\infty(\overline W)$ is
such that 
\begin{equation}\label{eq:2.3}
\nabla\psi\neq 0\quad\mbox{in}\quad \overline W
\end{equation}
 and $\lambda>0$ is a constant big enough. If $\Gamma\neq\emptyset$ we also suppose that
\begin{equation}\label{eq:2.4}
\frac{\partial\varphi}{\partial x_1}(0,x')>0,\quad\forall x'.
\end{equation}
If $\varphi=e^{\lambda\psi}$ the condition (\ref{eq:2.4}) is equivalent to
\begin{equation}\label{eq:2.5}
\frac{\partial\psi}{\partial x_1}(0,x')>0,\quad\forall x'.
\end{equation}
Let $\phi\in C^\infty(\overline W)$, supp$\,\phi\subset \overline W$, and let $u$ be as in Theorem 2.1. 
The next proposition follows from Propositions 1 and 2 of \cite{kn:LR1}.

\begin{prop} Let $\varphi$ satisfy (\ref{eq:2.2}). If $\Gamma\neq\emptyset$ we also suppose that $\varphi$ satisfies (\ref{eq:2.4}).
Then there exist constants $C, \hbar_0>0$ such that for all $0<\hbar\le\hbar_0$ we have the estimate
\begin{equation}\label{eq:2.6}
\int_X\left(|\phi u|^2+|\hbar\nabla(\phi u)|^2\right)e^{2\varphi/\hbar}dx\le C\hbar^3\int_X|\Delta_g(\phi u)|^2e^{2\varphi/\hbar}dx.
\end{equation}
\end{prop}

We take now $\hbar=\kappa h^{4/3}$, where $\kappa>0$ is a small parameter independent of $h$. By (\ref{eq:2.6}) we have
$$\int_X\left(|\phi u|^2+\kappa^2h^{2/3}|h\nabla(\phi u)|^2\right)e^{2\varphi/\kappa h^{4/3}}dx\le C\kappa^3\int_X|h^2\Delta_g(\phi u)|^2e^{2\varphi/\kappa h^{4/3}}dx$$
 $$\le C\kappa^3\int_X|(P(h)-z)(\phi u)|^2e^{2\varphi/\kappa h^{4/3}}dx+C\kappa^3\int_X|(V-z)(\phi u)|^2e^{2\varphi/\kappa h^{4/3}}dx$$
 $$\le C\kappa^3\int_X|(P(h)-z)(\phi u)|^2e^{2\varphi/\kappa h^{4/3}}dx+C(\sup|V|+C_0)^2\kappa^3\int_X|\phi u|^2e^{2\varphi/\kappa h^{4/3}}dx.$$
Taking $\kappa$ small enough we can absorb the last term in the right-hand side of the above inequality. Thus we obtain the following

\begin{prop} Let $\varphi$ satisfy (\ref{eq:2.2}). If $\Gamma\neq\emptyset$ we also suppose that $\varphi$ satisfies (\ref{eq:2.4}).
Then there exist constants $C, \kappa_0>0$ such that for all $0<\kappa\le\kappa_0$ and all $0<h\le 1$ we have the estimate
\begin{equation}\label{eq:2.7}
\int_X\left(|\phi u|^2+|h\nabla(\phi u)|^2\right)e^{2\varphi/\kappa h^{4/3}}dx\le C\kappa h^{-2/3}\int_X|(P(h)-z)(\phi u)|^2e^{2\varphi/\kappa h^{4/3}}dx.
\end{equation}
\end{prop}

In what follows in this section we will derive the estimate (\ref{eq:2.1}) from (\ref{eq:2.7}). Given a small parameter $\epsilon>0$,
independent of $h$, we denote $X_\epsilon=\{x\in X:{\rm dist}_g(x,\partial X)>\epsilon\}$ if $\partial X\neq \emptyset$, 
$X_\epsilon=X$ if $\partial X=\emptyset$. Taking $\epsilon$ small enough we can arrange that $U\subset X_\epsilon$. 
We will first derive from (\ref{eq:2.7}) the following

\begin{lemma} If $\partial X\neq \emptyset$, there exists a positive constant $\widetilde\gamma$ independent of $h$ such that 
for all $0<h\le 1$ we have the estimate
\begin{equation}\label{eq:2.8}
\|u\|_{H_h^1(X)}\le e^{\widetilde\gamma h^{-4/3}}\|(P(h)-z)u\|_{L^2(X)}+e^{\widetilde\gamma h^{-4/3}}\|u\|_{H_h^1(X_\epsilon)}.
\end{equation}
\end{lemma}

{\it Proof.} Let $\zeta\in C^\infty(\overline X)$ be such that $\zeta=1$ in $X\setminus X_\epsilon$, $\zeta =0$ in $X_{2\epsilon}$. 
Set $\psi(x)={\rm dist}_g(x,\partial X)$. Clearly, $\psi$ is $C^\infty$ smooth on supp$\,\zeta$, provided $\epsilon$ is small enough. 
Moreover, the function $\psi$ satisfies the conditions (\ref{eq:2.3}) and (\ref{eq:2.5}) on supp$\,\zeta$. Indeed, in the local coordinates $(x_1,x')$ above,
we have $\psi=x_1$. 
Let also $\eta_j\in C_0^\infty(\partial X)$, $j=1,...,J,$ be a partition of the unity on $\partial X$ such that the estimate 
(\ref{eq:2.7}) holds with $\varphi=e^{\lambda\psi}$, $\lambda\gg 1$, and 
$\phi$ replaced by $\phi_j=\zeta\eta_j$.
Taking into account that
$$[P(h),\phi_j]=-h^2[\Delta_g,\zeta\eta_j]=-h^2[\Delta_g,\eta_j]\zeta-h^2\eta_j[\Delta_g,\zeta]$$
and that $[\Delta_g,\zeta]$ is supported in $X_\epsilon\setminus X_{2\epsilon}$, 
we get from (\ref{eq:2.7})
$$\int_X\left(|\phi_j u|^2+|h\nabla(\phi_j u)|^2\right)e^{2\varphi/\kappa h^{4/3}}dx\le C\kappa h^{-2/3}\int_X|\zeta(P(h)-z)u|^2e^{2\varphi/\kappa h^{4/3}}dx$$
$$+C\kappa \int_X\left(|\zeta u|^2+|h\nabla(\zeta u)|^2\right)e^{2\varphi/\kappa h^{4/3}}dx+C\kappa \int_{X_\epsilon\setminus X_{2\epsilon}}\left(|u|^2+|h\nabla u|^2\right)e^{2\varphi/\kappa h^{4/3}}dx.$$
Summing up the above inequalities and using that $\zeta=\sum_{j=1}^J\phi_j$, we obtain
$$\int_X\left(|\zeta u|^2+|h\nabla(\zeta u)|^2\right)e^{2\varphi/\kappa h^{4/3}}dx\le C\kappa h^{-2/3}\int_X|\zeta(P(h)-z)u|^2e^{2\varphi/\kappa h^{4/3}}dx$$
$$+C\kappa \int_X\left(|\zeta u|^2+|h\nabla(\zeta u)|^2\right)e^{2\varphi/\kappa h^{4/3}}dx+C\kappa \int_{X_\epsilon\setminus X_{2\epsilon}}\left(|u|^2+|h\nabla u|^2\right)e^{2\varphi/\kappa h^{4/3}}dx$$
with a new constant $C>0$. 
Taking $\kappa$ small enough we can absorb the second term in the right-hand side of the above inequality and obtain the estimate
$$\int_X\left(|\zeta u|^2+|h\nabla(\zeta u)|^2\right)e^{2\varphi/\kappa h^{4/3}}dx\le 2C\kappa h^{-2/3}\int_X|\zeta(P(h)-z)u|^2e^{2\varphi/\kappa h^{4/3}}dx$$
$$+2C\kappa \int_{X_\epsilon\setminus X_{2\epsilon}}\left(|u|^2+|h\nabla u|^2\right)e^{2\varphi/\kappa h^{4/3}}dx.$$
Clearly, this implies
\begin{equation}\label{eq:2.9}
\|\zeta u\|_{H_h^1(X)}\le e^{\widetilde\gamma h^{-4/3}}\|(P(h)-z)u\|_{L^2(X)}+e^{\widetilde\gamma h^{-4/3}}\|u\|_{H_h^1(X_\epsilon\setminus X_{2\epsilon})}.
\end{equation}
with some constant $\widetilde\gamma>0$. Since
$$\|(1-\zeta) u\|_{H_h^1(X)}\lesssim \|u\|_{H_h^1(X_\epsilon)},$$
we get (\ref{eq:2.8}) from (\ref{eq:2.9}).
\eproof

Theorem 2.1 is a consequence of Lemma 2.4 and the following

\begin{lemma} Given any $\beta>0$ independent of $h$ there exists a positive constant $\gamma$ independent of $h$ such that 
for all $0<h\le 1$ we have the estimate
\begin{equation}\label{eq:2.10}
\|u\|_{H_h^1(X_\epsilon)}\le e^{\gamma h^{-4/3}}\|(P(h)-z)u\|_{L^2(X)}+e^{-\beta h^{-4/3}}\|u\|_{H_h^1(X)}
+e^{\gamma h^{-4/3}}\|u\|_{H_h^1(U)}.
\end{equation}
\end{lemma}

{\it Proof.} Given any $0<\rho\le\epsilon/6$ there are an integer $I=I(\rho)\ge 1$ and balls 
$B_i(\rho)=\{x\in X: {\rm dist}_g(x,y_i)<\rho\}$,
$i=1,...,I$, $y_1\in U$, $y_i\in X_\epsilon$, $i=2,...,I$, such that $X_\epsilon\subset\cup_{i=1}^IB_i(\rho)$. 
If $\partial X\neq \emptyset$, clearly $\partial X\cap\overline{B_i(5\rho)}=\emptyset$, $i=1,...,I$. 
Taking $\rho$ small enough we can
also arrange that $B_1(\rho)\subset U$. Set $\psi(x)=-{\rm dist}_g(x,y_i)\in C^\infty(B_i(5\rho)\setminus\{y_i\})$ 
and let
$\phi\in C_0^\infty(B_i(5\rho)\setminus B_i(\rho/2))$ be such that $\phi=1$ in $B_i(4\rho)\setminus B_i(\rho)$. 
Clearly, the function
$\psi$ is smooth on supp$\,\phi$ and satisfies the condition (\ref{eq:2.3}). Thus, since 
supp$\,\phi\cap \partial X=\emptyset$, we can apply
the estimate (\ref{eq:2.7}) with $\varphi=e^{\lambda\psi}$, $\lambda\gg 1$, to obtain
$$e^{2e^{-3\lambda\rho}/\kappa h^{4/3}}\int_{B_i(3\rho)\setminus B_i(\rho)}\left(|u|^2+|h\nabla u|^2\right)dx$$
$$\le\int_{B_i(3\rho)\setminus B_i(\rho)}\left(|u|^2+|h\nabla u|^2\right)e^{2\varphi/\kappa h^{4/3}}dx$$
$$\le\int_X\left(|\phi u|^2+|h\nabla(\phi u)|^2\right)e^{2\varphi/\kappa h^{4/3}}dx$$
$$\lesssim h^{-2/3}\int_X|(P(h)-z)(\phi u)|^2e^{2\varphi/\kappa h^{4/3}}dx$$
$$\lesssim h^{-2/3}\int_X|\phi(P(h)-z)u|^2e^{2\varphi/\kappa h^{4/3}}dx$$
$$+\int_{(B_i(5\rho)\setminus B_i(4\rho))\cup
(B_i(\rho)\setminus B_i(\rho/2))}\left(|u|^2+|h\nabla u|^2\right)e^{2\varphi/\kappa h^{4/3}}dx$$
$$\lesssim h^{-2/3}e^{2e^{-\lambda\rho/2}/\kappa h^{4/3}}\int_X|(P(h)-z)u|^2dx$$
$$+ e^{2e^{-4\lambda\rho}/\kappa h^{4/3}}\int_{B_i(5\rho)\setminus B_i(4\rho)}\left(|u|^2+|h\nabla u|^2\right)dx$$
$$+ e^{2e^{-\lambda\rho/2}/\kappa h^{4/3}}\int_{B_i(\rho)\setminus B_i(\rho/2)}\left(|u|^2+|h\nabla u|^2\right)dx.$$
This implies 
\begin{equation}\label{eq:2.11}
\|u\|_{H_h^1(B_i(3\rho))}^2\lesssim h^{-2/3}e^{2c_1\kappa^{-1}h^{-4/3}}\|(P(h)-z)u\|_{L^2(X)}^2$$
$$+e^{-2c_2\kappa^{-1}h^{-4/3}}\|u\|_{H_h^1(X)}^2+e^{2c_1\kappa^{-1} h^{-4/3}}\|u\|_{H_h^1(B_i(\rho))}^2
\end{equation}
for every $0<\kappa\ll 1$ independent of $h$, where $c_1=e^{-\lambda\rho/2}-e^{-3\lambda\rho}>0$ and 
$c_2=e^{-3\lambda\rho}-e^{-4\lambda\rho}>0$. 
Choosing the parameter $\kappa$ suitably we will show now that (\ref{eq:2.11}) implies the estimate
\begin{equation}\label{eq:2.12}
\|u\|_{H_h^1(B_i(\rho))}^2\lesssim e^{2\gamma_ih^{-4/3}}\|(P(h)-z)u\|_{L^2(X)}^2$$
$$+e^{-2\beta h^{-4/3}}\|u\|_{H_h^1(X)}^2+e^{2\gamma_i h^{-4/3}}\|u\|_{H_h^1(B_1(\rho))}^2
\end{equation}
for all $i=1,...,I,$ and for any $\beta>0$ independent of $h$ with some constant $\gamma_i>0$ depending on $\beta$.
The estimate (\ref{eq:2.12}) is trivial for $i=1$. Let $i\ge 2$. Since $X$ is connected, there exist integers
$i_1,...,i_L$, $2\le L\le I$, $i_1=1$, $i_L=i$, $2\le i_\ell\le I$ if $2\le\ell\le L-1$ such that
\begin{equation}\label{eq:2.13}
B_{i_{\ell-1}}(\rho)\cap B_{i_{\ell}}(\rho)\neq\emptyset,\quad 2\le\ell\le L.
\end{equation}
Clearly, (\ref{eq:2.13}) implies
\begin{equation}\label{eq:2.14}
B_{i_{\ell}}(\rho)\subset B_{i_{\ell-1}}(3\rho),\quad 2\le\ell\le L.
\end{equation}
We now apply the estimate (\ref{eq:2.11}) with $i$ replaced by $i_{\ell-1}$ and $\kappa$ replaced by 
$\kappa_\ell$ to be chosen later on. Thus, in view of (\ref{eq:2.14}), we get
\begin{equation}\label{eq:2.15}
\|u\|_{H_h^1(B_{i_\ell}(\rho))}^2\lesssim h^{-2/3}e^{2c_1\kappa_\ell^{-1}h^{-4/3}}\|(P(h)-z)u\|_{L^2(X)}^2$$
$$+e^{-2c_2\kappa_\ell^{-1}h^{-4/3}}\|u\|_{H_h^1(X)}^2
+e^{2c_1\kappa_\ell^{-1} h^{-4/3}}\|u\|_{H_h^1(B_{i_{\ell-1}}(\rho))}^2
\end{equation}
for all $\ell=2,...,L$. Iterating these inequalities leads to the estimate
\begin{equation}\label{eq:2.16}
\|u\|_{H_h^1(B_{i_L}(\rho))}^2\lesssim h^{-2/3}Q_1\|(P(h)-z)u\|_{L^2(X)}^2+Q_2\|u\|_{H_h^1(X)}^2
+Q_3\|u\|_{H_h^1(B_{i_1}(\rho))}^2
\end{equation}
where 
$$Q_1=\sum_{\ell=2}^L\exp\left(2h^{-4/3}\sum_{\nu=\ell}^L\frac{c_1}{\kappa_\nu}\right),$$
$$Q_2=\exp\left(-2h^{-4/3}\frac{c_2}{\kappa_L}\right),$$
if $L=2$,
$$Q_2=\exp\left(-2h^{-4/3}\frac{c_2}{\kappa_L}\right)+\sum_{\ell=2}^{L-1}\exp\left(-2h^{-4/3}\frac{c_2}{\kappa_\ell}+2h^{-4/3}\sum_{\nu=\ell+1}^L\frac{c_1}{\kappa_\nu}
\right),$$
if $L\ge 3$, and 
$$Q_3=\exp\left(2h^{-4/3}\sum_{\nu=2}^L\frac{c_1}{\kappa_\nu}\right).$$
Observe now that given any $\beta>0$ we can choose the parameters $\kappa_\ell$, $\ell=2,...,L,$ small enough in order to arrange
the inequalities
$$\frac{c_2}{\kappa_L}\ge\beta,\quad \frac{c_2}{\kappa_\ell}\ge \beta+\sum_{\nu=\ell+1}^L\frac{c_1}{\kappa_\nu}$$
for every $2\le\ell\le L-1$ (if $L\ge 3$). Therefore the estimate (\ref{eq:2.12}) follows from (\ref{eq:2.16}).
Finally, observe that summing up all the inequalities (\ref{eq:2.12}) leads to the estimate (\ref{eq:2.10}) with any 
$$\gamma>\max_{1\le i\le I}\gamma_i.$$ 
\eproof

Combining the estimates (\ref{eq:2.8}) and (\ref{eq:2.10}) we get
$$\|u\|_{H_h^1(X)}\le 2e^{(\widetilde\gamma+\gamma)h^{-4/3}}\|(P(h)-z)u\|_{L^2(X)}$$
$$+e^{-(\beta-\widetilde\gamma) h^{-4/3}}\|u\|_{H_h^1(X)}
+e^{(\widetilde\gamma+\gamma)h^{-4/3}}\|u\|_{H_h^1(U)}.$$
Clearly, taking $\beta$ big enough we can absorb the second term in the right-hand side of the above inequality and obtain 
(\ref{eq:2.1}) with a new constant $\gamma$.
\eproof

\section{Construction of the phase and weight functions on $Y$} 

We will first construct the weight function. In what follows $b>0$ will be a parameter independent of $h$ to be fixed
in the proof of Lemma 4.2 depending only on the dimension $n$, the Riemannian metric $\omega$ and the constants $C$
appearing in the conditions (\ref{eq:1.5}) and (\ref{eq:1.6}). Since the function $f$ is increasing, there is $r_1\ge r_0$
depending on $b$ such that $f(r)\ge 2b$ for all $r\ge r_1$. If $V$ is of compact support we take $r_1$ large enough to assure that
$V=0$ in $Y_{r_1}$. With this in mind we introduce the continuous function
$$\mu(r)=
\left\{
\begin{array}{lll}
 (f(r)-b)^2&\mbox{for}& r_1\le r\le a,\\
 (f(a)-b)^2+a^{-2s+1}-r^{-2s+1}&\mbox{for}& r\ge a,
\end{array}
\right.
$$
where
\begin{equation}\label{eq:3.1}
s=\frac{1+\epsilon}{2},\quad \epsilon=\left(\log\frac{1}{h}\right)^{-1},
\end{equation}
and $a=h^{-m}$ with
$$m=m_0+\frac{\epsilon(\lambda+m_0+t)}{\delta-1},$$
$m_0$ being as in Section 1, $\lambda=\log\log\frac{1}{h}$. If $V$ is of compact support we set 
$$m=m_0+\epsilon t.$$
Here $t>0$  
is a parameter independent of $h$ to be fixed in the proof of Lemma 3.3.
Clearly, the first derivative (in sense of distributions) of $\mu$ satisfies
$$\mu'(r)=
\left\{
\begin{array}{lll}
 2f'(r)(f(r)-b)&\mbox{for}& r_1\le r<a,\\
 (2s-1)r^{-2s}&\mbox{for}& r>a.
\end{array}
\right.
$$

\begin{lemma} For all $r\ge r_1$, $r\neq a$, we have the bounds
\begin{equation}\label{eq:3.2}
\frac{\mu(r)}{\mu'(r)}\lesssim  \epsilon^{-1}f(a)^2r^{2s}, 
\end{equation}
\begin{equation}\label{eq:3.3}
\frac{\mu(r)^2}{\mu'(r)}\lesssim  \epsilon^{-1}f(a)^4r^{2s}. 
\end{equation}
\end{lemma}

{\it Proof.} For $r_1\le r<a$ we have the bounds
\begin{equation}\label{eq:3.4}
\frac{\mu(r)}{\mu'(r)}=\frac{f(r)-b}{2f'(r)}<\frac{f(r)}{2f'(r)}\lesssim r, 
\end{equation}
\begin{equation}\label{eq:3.5}
\frac{\mu(r)^2}{\mu'(r)}=\frac{(f(r)-b)^2}{2f'(r)}<\frac{f(r)^2}{2f'(r)}\lesssim f(a)r. 
\end{equation}
For $r>a$ we have $\mu={\cal O}(f(a)^2)$ and $\mu'(r)=\epsilon r^{-2s}$. 
\eproof

We now turn to the 
construction of the phase function
$\varphi\in C^1([r_1,+\infty))$ such that $\varphi(0)=0$ and $\varphi(r)>0$ for $r>0$. 
We define the first derivative of $\varphi$ by
$$\varphi'(r)=
\left\{
\begin{array}{lll}
 \tau f(r)^{-1}- \tau f(a)^{-1}&\mbox{for}& r_1\le r\le a,\\
 0&\mbox{for}& r\ge a,
\end{array}
\right.
$$
where 
\begin{equation}\label{eq:3.6}
\tau=\tau_0h^{-1/3}
\end{equation}
with a parameter $\tau_0\gg 1$ independent of $h$.
Clearly, the first derivative of $\varphi'$ satisfies
$$\varphi''(r)=
\left\{
\begin{array}{lll}
 -\tau f'(r)f(r)^{-2}&\mbox{for}& r_1\le r<a,\\
 0&\mbox{for}& r>a.
\end{array}
\right.
$$

\begin{lemma} If $f$ is given by (\ref{eq:1.3}) with $k<1$ we have the bounds
\begin{equation}\label{eq:3.7}
h^{-1}\varphi(r)\lesssim 
\left\{
\begin{array}{lll}
 h^{-4/3-m_0(1-k)}\left(\log(h^{-1})\right)^{\frac{1-k}{\delta-1}}&\mbox{if $V$ satisfies 
 (\ref{eq:1.9})},\\
 h^{-\frac{2(k+1)}{3k}}&\mbox{if $V$ is of compact support},\\
\end{array}
\right.
\end{equation}
for all $r\ge r_1$. In the other two cases we have the bounds
\begin{equation}\label{eq:3.8}
h^{-1}\varphi(r)\lesssim 
\left\{
\begin{array}{lll}
 h^{-4/3}\log(h^{-1})&\mbox{if $f$ is given by
 (\ref{eq:1.3}) with $k=1$},\\
 h^{-4/3}&\mbox{if $f$ is given by (\ref{eq:1.3}) with $k>1$ or by 
 (\ref{eq:1.4})}.
\end{array}
\right.
\end{equation}
\end{lemma}

{\it Proof.} We have
$$\max_{r\ge r_1}\varphi=\int_{r_1}^a\varphi'(r)dr\le \tau\int_{r_1}^a \frac{dr}{f(r)}$$
$$\lesssim  \left\{
\begin{array}{lll}
 h^{-1/3-m(1-k)}&\mbox{if $f$ is given by
 (\ref{eq:1.3}) with $k<1$},\\
 h^{-1/3}\log(h^{-1})&\mbox{if $f$ is given by
 (\ref{eq:1.3}) with $k=1$},\\
 h^{-1/3}&\mbox{if $f$ is given by (\ref{eq:1.3}) with $k>1$ or by 
 (\ref{eq:1.4})}.
\end{array}
\right.
$$
Observe now that in the first case we have
$$m(1-k)=m_0(1-k)+\frac{1-k}{\delta-1}\epsilon\lambda +{\cal O}(\epsilon),$$
while if $V$ is of compact support we have
$$m(1-k)=m_0(1-k)+{\cal O}(\epsilon)=\frac{2(1-k)}{3k}+{\cal O}(\epsilon).$$
This clearly implies (\ref{eq:3.7}).
\eproof

For $r\ge r_1$, $r\neq a$, set
$$A(r)=\left(\mu\varphi'^2\right)'(r)$$
and
$$B(r)=\frac{\left(\mu(r)\left(h^{-1}r^{-\delta}f(r)^{-2}+hr^{-1}f(r)^{-2}+|\varphi''(r)|\right)\right)^2}{h^{-1}\varphi'(r)\mu(r)+\mu'(r)}.$$
If $V$ is of compact support we set
$$B(r)=\frac{\left(\mu(r)\left(hr^{-1}f(r)^{-2}+|\varphi''(r)|\right)\right)^2}{h^{-1}\varphi'(r)\mu(r)+\mu'(r)}.$$
The following lemma will play a crucial role in the proof of the Carleman estimates in the next section.

\begin{lemma} Given any $C>0$ independent of the variable $r$ and the parameters $h$, $\tau$ and $a$, there exist $\tau_1=\tau_1(C)>0$ and 
$h_0=h_0(C)>0$ so that for $\tau$ satisfying (\ref{eq:3.6}) with $\tau_0\ge\tau_1$ and for all $0<h\le h_0$ we have the inequality
\begin{equation}\label{eq:3.9}
A(r)-CB(r)-h^2(\mu q_0)'(r)\ge -\frac{E}{2}\mu'(r)
\end{equation}
for all $r\ge r_1$, $r\neq a$.
\end{lemma}

{\it Proof.} We will first bound from above the function
$$(\mu q_0)'=\left(q_0+\frac{\mu}{\mu'}q'_0\right)\mu'$$
using that $q_0$ satisfies the condition (\ref{eq:1.7}).
For $r_1\le r<a$ we have
$$q_0+\frac{\mu}{\mu'}q'_0\lesssim 1+r^{-1}f'(r)^{-1}f(r)^{-1}\lesssim 1.$$
For $r>a$, in view of (\ref{eq:3.2}), we have
$$q_0+\frac{\mu}{\mu'}q'_0\lesssim 1+\epsilon^{-1}f(a)^2r^\epsilon f(r)^{-2}.$$
Observe now that for $\epsilon$ small enough the function $r^\epsilon f(r)^{-2}$ is decreasing. Hence
$$r^\epsilon f(r)^{-2}\le a^\epsilon f(a)^{-2}\lesssim f(a)^{-2}$$
where we have used that $a^\epsilon={\cal O}(1)$. Thus we get the inequality
\begin{equation}\label{eq:3.10}
h^2(\mu q_0)'(r)\lesssim h\mu'(r)\le \frac{E}{8}\mu'(r)
\end{equation}
provided $h$ is small enough.

We will now bound from below the function $A(r)$ for $r_1\le r<a$. We have
$$A(r)=2\tau\varphi'(r)(f(r)-b)\partial_r\left(1-bf(r)^{-1}-(f(r)-b)f(a)^{-1}\right)$$ 
$$=2\tau\varphi'(r)(f(r)-b)\left(bf'(r)f(r)^{-2}-f'(r)f(a)^{-1}\right)$$
$$\ge b\tau\varphi'(r)f'(r)f(r)^{-1}-2\tau\varphi'(r)f'(r)(f(r)-b)f(a)^{-1}$$
$$\ge b\tau\varphi'(r)f'(r)f(r)^{-1}-\tau^2 f(r_1)^{-1}f(a)^{-1}\mu'(r).$$
Observe now that when $f$ is given by (\ref{eq:1.3}) we have 
$$\tau^2f(a)^{-1}=\tau^2 a^{-k}\lesssim h^{mk-2/3}$$
$$\lesssim 
\left\{
\begin{array}{lll}
 \epsilon^{k/(\delta-1)}&\mbox{if $V$ satisfies 
 (\ref{eq:1.9})},\\
 e^{-kt}&\mbox{if $V$ is of compact support},\\
\end{array}
\right.
$$
while when $f$ is given by (\ref{eq:1.4}) we have 
$$\tau^2f(a)^{-1}=\tau^2e^{-a^\alpha}\lesssim h^{-2/3}e^{-h^{-m\alpha}}\lesssim h.$$
Thus, taking $h$ small enough and $t$ big enough, we can arrange that the inequality
\begin{equation}\label{eq:3.11}
A(r)\ge b\tau\varphi'(r)f'(r)f(r)^{-1}-\frac{E}{8}\mu'(r)
\end{equation}
holds for all $r_1\le r<a$. 

We will now bound from above the function $B$ in the general case. When $V$ is of compact support the analysis of $B$
is much easier and we omit the details. 

Let first $r_1\le r\le \frac{a}{2}$. In this case we have
$$\varphi'(r)\ge C\tau f(r)^{-1}$$
with some constant $C>0$. 
Thus we obtain
$$B(r)\lesssim \frac{\mu(r)\left(h^{-2}r^{-2\delta}f(r)^{-4}+h^2r^{-2}f(r)^{-4}+\varphi''(r)^2\right)}{h^{-1}\varphi'(r)}$$ 
 $$\lesssim (\tau h)^{-1}\frac{\mu(r)r^{-2\delta}f'(r)^{-1}f(r)^{-3}}{\varphi'(r)^2}\tau\varphi'(r)f'(r)f(r)^{-1}$$
 $$+ h^3\frac{\mu(r)r^{-2}f(r)^{-4}}{\mu'(r)\varphi'(r)}\mu'(r)
 + h\frac{\mu(r)\varphi''(r)^2}{\mu'(r)\varphi'(r)}\mu'(r)$$ 
$$\lesssim \tau^{-3} h^{-1}r^{-2\delta}f'(r)^{-1}f(r)\tau\varphi'(r)f'(r)f(r)^{-1}$$
$$+ h^3\tau^{-1}f(r)^{-3}\mu'(r)
 + h\tau f'(r)f(r)^{-2}\mu'(r)$$ 
$$\lesssim \tau_0^{-3}\tau\varphi'(r)f'(r)f(r)^{-1} + h^{2/3}\mu'(r)$$
where we have used that $f'={\cal O}(f)$, $f^{-1}={\cal O}(1)$ together with the bound (\ref{eq:3.4}). 
The above bound together with (\ref{eq:3.10}) and (\ref{eq:3.11}) clearly imply (\ref{eq:3.9}), provided $\tau_0^{-1}$ and $h$ are taken small enough
depending on $C$.

Let now $\frac{a}{2}<r<a$. In view of (\ref{eq:3.4}), we have 
$$B(r)\le\left(\frac{\mu(r)}{\mu'(r)}\right)^2\left(h^{-1}r^{-\delta}f(r)^{-2}+hr^{-1}f(r)^{-2}+|\varphi''(r)|\right)^2\mu'(r)$$
$$\lesssim \left(h^{-1}r^{1-\delta}f(r)^{-2}+hf(r)^{-2}+\tau f(r)^{-1}\right)^2\mu'(r)$$
$$\lesssim \left(h^{-1}f(a/2)^{-2}+\tau f(a/2)^{-1}\right)^2\mu'(r)\lesssim h^{2/3}\mu'(r).$$ 
 Again, this bound together with (\ref{eq:3.10}) and (\ref{eq:3.11}) imply (\ref{eq:3.9}).

It remains to consider the case $r>a$. Taking into account
that $s$ satisfies (\ref{eq:3.1}) and using the bound (\ref{eq:3.2}), we get
$$B(r)=\left(\frac{\mu(r)}{\mu'(r)}\right)^2\left(h^{-2}r^{-2\delta}f(r)^{-4}+h^{2}r^{-2}f(r)^{-4}\right)\mu'(r)$$
$$\lesssim \epsilon^{-2}\left(h^{-2}f(a)^{4}r^{-2\delta+4s}f(r)^{-4}+h^2f(a)^{4}r^{-2+4s}f(r)^{-4}\right)\mu'(r)$$
$$\lesssim \epsilon^{-2}\left(h^{-2}a^{-2\delta+4s}+h^2a^{-2+4s}\right)\mu'(r)$$
$$\lesssim \epsilon^{-2}\left(h^{2m(\delta-2s)-2}+h^{-2m\epsilon+2}\right)\mu'(r)$$
$$\lesssim \left(h^{2m(\delta-1-\epsilon)-2-2\epsilon\lambda}+\epsilon^{-2}h^{2}\right)\mu'(r).$$
On the other hand, we have
$$m(\delta-1-\epsilon)-1-\epsilon\lambda=\left(m_0+\frac{\epsilon(\lambda+m_0+t)}{\delta-1}\right)(\delta-1-\epsilon)-1-\epsilon\lambda$$ 
$$=(\delta-1)m_0-1+\epsilon t-{\cal O}(\lambda\epsilon^2)\ge \epsilon t/2.$$
Hence 
\begin{equation}\label{eq:3.12}
B(r)\lesssim\left(e^{-t}+h\right)\mu'(r)\le \frac{E}{4C}\mu'(r)
\end{equation}
provided $h$ is taken small enough and $t$ big enough, independent of $h$.
 Since in this case $A(r)=0$, the bound (\ref{eq:3.12}) together with (\ref{eq:3.10}) clearly imply (\ref{eq:3.9}).
\eproof

\section{Carleman estimates on $Y_{r_1}$} 

Our goal in this section is to prove the following

\begin{Theorem} Let $s$ satisfy (\ref{eq:3.1}). Then, under the conditions of Theorem 1.1, for all functions
$u\in H^2(Y_{r_1},d{\rm Vol}_g)$ such that $r^s(P(h)-E\pm i\varepsilon)u\in L^2(Y_{r_1},d{\rm Vol}_g)$, $u=\partial_ru=0$ on $\partial Y_{r_1}$, and for all
$0<h\ll 1$, we have the estimate 
 \begin{equation}\label{eq:4.1}
\|r^{-s}e^{\varphi/h}u\|_{L^2(Y_{r_1},d{\rm Vol}_g)}+\|r^{-s}e^{\varphi/h}{\cal D}_ru\|_{L^2(Y_{r_1},d{\rm Vol}_g)}$$
$$\le Cf(a)^2(\epsilon h)^{-1}\|r^se^{\varphi/h}(P(h)-E\pm i\varepsilon)u\|_{L^2(Y_{r_1},d{\rm Vol}_g)}$$ 
$$+C\tau f(a)\varepsilon^{1/2}(\epsilon h)^{-1/2}\|e^{\varphi/h}u\|_{L^2(Y_{r_1},d{\rm Vol}_g)}
\end{equation}
with a constant $C>0$ independent of $h$, $\varepsilon$ and $u$, where ${\cal D}_r:=-ih\partial_r$.
\end{Theorem}

{\it Proof.} In what follows we denote by $\|\cdot\|$ and $\langle\cdot,\cdot\rangle$
the norm and the scalar product in $L^2(S)$. Note that $d{\rm Vol}_g=p(r,\theta)drd\theta$ on $Y_{r_1}$. 
Set $v=p^{1/2}e^{\varphi/h}u$ and
$${\cal P}^\pm(h)=p^{1/2}(P(h)-E\pm i\varepsilon)p^{-1/2},$$
$${\cal P}^\pm_\varphi(h)=e^{\varphi/h}{\cal P}^\pm(h)e^{-\varphi/h}.$$
Using (\ref{eq:1.1}) we can write the operator ${\cal P}^\pm(h)$ as follows
$${\cal P}^\pm(h)={\cal D}_r^2+L_\theta(r)-E\pm i\varepsilon +V+h^2q$$
where we have put $L_\theta(r)=-h^2\Lambda_\theta(r)\ge 0$. Since the function $\varphi$
depends only on the variable $r$, this implies
$${\cal P}^\pm_\varphi(h)={\cal D}_r^2+L_\theta(r)-E\pm i\varepsilon -\varphi'^2+h\varphi''+
2i\varphi'{\cal D}_r+V+h^2q.$$
For $r\ge r_1$, $r\neq a$, introduce the function
$$F(r)=-\langle (L_\theta(r)-E-\varphi'(r)^2+h^2q_0)v(r,\cdot),v(r,\cdot)\rangle+\|{\cal D}_rv(r,\cdot)\|^2$$
and observe that its first derivative is given by
$$F'(r)=-\langle [\partial_r,L_\theta(r)]v(r,\cdot),v(r,\cdot)\rangle
+((\varphi')^2-h^2q_0)'\|v(r,\cdot)\|^2$$
$$-2h^{-1}{\rm Im}\,\langle {\cal P}^\pm_\varphi(h)v(r,\cdot),{\cal D}_rv(r,\cdot)\rangle$$
$$\pm 2\varepsilon h^{-1}{\rm Re}\,\langle v(r,\cdot),{\cal D}_rv(r,\cdot)\rangle+4h^{-1}\varphi'\|{\cal D}_rv(r,\cdot)\|^2$$ 
$$+2h^{-1}{\rm Im}\,\langle (V+h\varphi''+h^2(q-q_0))v(r,\cdot),{\cal D}_rv(r,\cdot)\rangle.$$
Thus, if $\mu$ is the weight function defined in the previous section, we obtain the identity
$$\mu'F+\mu F'=-\langle(\mu[\partial_r,L_\theta(r)]+\mu'L_\theta(r))v(r,\cdot),v(r,\cdot)\rangle$$ 
$$+(E\mu'+(\mu(\varphi')^2-h^2\mu q_0)')\|v(r,\cdot)\|^2$$
$$-2h^{-1}\mu{\rm Im}\,\langle {\cal P}^\pm_\varphi(h)v(r,\cdot),{\cal D}_rv(r,\cdot)\rangle$$
$$\pm 2\varepsilon h^{-1}\mu{\rm Re}\,\langle v(r,\cdot),{\cal D}_rv(r,\cdot)\rangle+(\mu'+4h^{-1}\varphi'\mu)\|{\cal D}_rv(r,\cdot)\|^2$$ 
$$+2h^{-1}\mu{\rm Im}\,\langle (V+h\varphi''+h^2(q-q_0))v(r,\cdot),{\cal D}_rv(r,\cdot)\rangle.$$
We need now the following

\begin{lemma} For all $r\ge r_1$, $r\neq a$, we have the inequality
\begin{equation}\label{eq:4.2}
\langle(\mu[\partial_r,L_\theta(r)]+\mu'L_\theta(r))v,v\rangle\le 0,\quad\forall v\in H^1(S).
\end{equation}
\end{lemma}

{\it Proof.} Clearly, the operator in the left-hand side of (\ref{eq:4.2}) is of the form
$$-h^2\sum_{i,j}\partial_{\theta_i}(\Phi^{ij}(r,\theta)\partial_{\theta_j})$$
where
$$\Phi^{ij}=\mu\partial_rg^{ij}+\mu'g^{ij}=\mu\partial_r(g^{ij}-f^{-2}\omega^{ij})+\mu'(g^{ij}-f^{-2}\omega^{ij})
+(\mu f^{-2})'\omega^{ij}.$$
Thus the left-hand side of (\ref{eq:4.2}) can be written in the form
$$h^2\sum_{i,j}\langle\Phi^{ij}\partial_{\theta_i}v,\partial_{\theta_j}v\rangle.$$
Therefore, to prove (\ref{eq:4.2}) it suffices to show that
\begin{equation}\label{eq:4.3}
\sum_{i,j}\Phi^{ij}\xi_i\overline{\xi_j}\le 0,\quad\forall \xi\in \mathbb{C}^{n-1}.
\end{equation}
To this end, we will use the conditions (\ref{eq:1.5}) and (\ref{eq:1.6}). For $r_1\le r<a$ we have 
$$(\mu f^{-2})'=-2\left(1-\frac{b}{f(r)}\right)\frac{bf'(r)}{f(r)^2}\le -\frac{bf'(r)}{f(r)^2}.$$
Observe now that the function $f$ satisfies the inequality
\begin{equation}\label{eq:4.4}
f'(r)\ge\widetilde Cr^{-1}f(r),\quad\widetilde C>0.
\end{equation}
In view of (\ref{eq:4.4}), for $r>a$, we have
$$(\mu f^{-2})'=-\frac{2\mu(r)f'(r)}{f(r)^3}+\frac{\mu'(r)}{f(r)^2}$$
$$\le -\frac{\mu(r)f'(r)}{f(r)^3}-\frac{\widetilde C(f(a)-b)^2-2s+1}{rf(r)^2}\le -\frac{\mu(r)f'(r)}{f(r)^3}$$
provided $a$ is taken large enough. 
Thus, using that
$$\sum_{i,j}\omega^{ij}\xi_i\overline{\xi_j}\ge C_\sharp|\xi|^2,\quad C_\sharp>0,$$
 we obtain with some constant $C>0$ independent of the parameter $b$,
$$\sum_{i,j}\Phi^{ij}\xi_i\overline{\xi_j}\le -\frac{bf'(r)}{f(r)^2}\sum_{i,j}\omega^{ij}\xi_i\overline{\xi_j}+\frac{C\mu(r)f'(r)}{f(r)^4}|\xi|^2
+\frac{C\mu'(r)}{f(r)^3}|\xi|^2$$
$$\le -\frac{C_\sharp bf'(r)}{f(r)^2}|\xi|^2+\frac{3Cf'(r)}{f(r)^2}|\xi|^2=-\frac{Cf'(r)}{f(r)^2}|\xi|^2\le 0$$
for $r_1\le r<a$, if we choose $b=4C/C_\sharp$. In view of (\ref{eq:4.4}), for $r>a$ we have
$$\sum_{i,j}\Phi^{ij}\xi_i\overline{\xi_j}\le -\frac{C_\sharp\mu(r)f'(r)}{f(r)^3}|\xi|^2+\frac{C\mu(r)f'(r)}{f(r)^4}|\xi|^2
+\frac{C\mu'(r)}{f(r)^3}|\xi|^2$$
$$\le -\frac{C_\sharp\mu(r)f'(r)}{2f(r)^3}|\xi|^2\le 0,$$
provided $a$ is taken large enough. Thus in both cases we get (\ref{eq:4.3}).
\eproof

Using (\ref{eq:4.2}) we get the inequality
$$\mu'F+\mu F'\ge (E\mu'+(\mu(\varphi')^2-h^2\mu q_0)')\|v(r,\cdot)\|^2+(\mu'+4h^{-1}\varphi'\mu)\|{\cal D}_rv(r,\cdot)\|^2$$
$$-\frac{3h^{-2}\mu^2}{\mu'}\|{\cal P}^\pm_\varphi(h)v(r,\cdot)\|^2-\frac{\mu'}{3}\|{\cal D}_rv(r,\cdot)\|^2$$
$$-\varepsilon h^{-1}\mu\left(\|v(r,\cdot)\|^2+\|{\cal D}_rv(r,\cdot)\|^2\right)$$
$$-3h^{-2}\mu^2(\mu'+4h^{-1}\varphi'\mu)^{-1}\|(V+h\varphi''+h^2(q-q_0))v(r,\cdot)\|^2$$ 
$$-\frac{1}{3}(\mu'+4h^{-1}\varphi'\mu)\|{\cal D}_rv(r,\cdot)\|^2$$
 $$\ge \left(E\mu'+(\mu(\varphi')^2-h^2\mu q_0)'\right)\|v(r,\cdot)\|^2$$
 $$-C\mu^2(\mu'+h^{-1}\varphi'\mu)^{-1}(h^{-1}r^{-\delta}f(r)^{-2}+hr^{-\beta}f(r)^{-2}
 +|\varphi''|)^2\|v(r,\cdot)\|^2$$
 $$+\frac{\mu'}{3}\|{\cal D}_rv(r,\cdot)\|^2
-\frac{3h^{-2}\mu^2}{\mu'}\|{\cal P}^\pm_\varphi(h)v(r,\cdot)\|^2$$ 
$$-\varepsilon h^{-1}\mu\left(\|v(r,\cdot)\|^2+\|{\cal D}_rv(r,\cdot)\|^2\right)$$
with some constant $C>0$. Now we use Lemma 3.3 to conclude that
$$\mu'F+\mu F'\ge \frac{E}{2}\mu'\|v(r,\cdot)\|^2+\frac{\mu'}{3}\|{\cal D}_rv(r,\cdot)\|^2-\frac{3h^{-2}\mu^2}{\mu'}\|{\cal P}^\pm_\varphi(h)v(r,\cdot)\|^2$$ $$
-\varepsilon h^{-1}\mu\left(\|v(r,\cdot)\|^2+\|{\cal D}_rv(r,\cdot)\|^2\right).$$
We now integrate this inequality with respect to $r$. Since $F(r_1)=0$, we have
$$\int_{r_1}^\infty(\mu'F+\mu F')dr=0.$$
Thus we obtain the estimate
\begin{equation}\label{eq:4.5}
\frac{E}{2}\int_{r_1}^\infty\mu'\|v(r,\cdot)\|^2dr+\frac{1}{3}\int_{r_1}^\infty\mu'\|{\cal D}_rv(r,\cdot)\|^2dr$$ 
$$\le 3h^{-2}\int_{r_1}^\infty\frac{\mu^2}{\mu'}
\|{\cal P}^\pm_\varphi(h)v(r,\cdot)\|^2dr$$ $$
+\varepsilon h^{-1}\int_{r_1}^\infty\mu\left(\|v(r,\cdot)\|^2+\|{\cal D}_rv(r,\cdot)\|^2\right)dr.
\end{equation}
Using that $\mu={\cal O}(f(a)^2)$, $\mu'\ge \epsilon r^{-2s}$ together with (\ref{eq:3.3}) we get from (\ref{eq:4.5})
\begin{equation}\label{eq:4.6}
\epsilon\int_{r_1}^\infty r^{-2s}\left(\|v(r,\cdot)\|^2+\|{\cal D}_rv(r,\cdot)\|^2\right)dr\le Cf(a)^4h^{-2}\epsilon^{-1}\int_{r_1}^\infty r^{2s}\|{\cal P}^\pm_\varphi(h)v(r,\cdot)\|^2dr$$ $$
+C\varepsilon h^{-1}f(a)^2\int_{r_1}^\infty\left(\|v(r,\cdot)\|^2+\|{\cal D}_rv(r,\cdot)\|^2\right)dr
\end{equation}
with some constant $C>0$ independent of $h$ and $\varepsilon$. On the other hand, we have the identity
$${\rm Re}\,\int_{r_1}^\infty\langle 2i\varphi'{\cal D}_rv(r,\cdot),v(r,\cdot)\rangle dr=\int_{r_1}^\infty h\varphi''\|v(r,\cdot)\|^2dr$$
and hence
$${\rm Re}\,\int_{r_1}^\infty\langle {\cal P}^\pm_\varphi(h)v(r,\cdot),v(r,\cdot)\rangle dr
=\int_{r_1}^\infty\|{\cal D}_rv(r,\cdot)\|^2dr
+\int_{r_1}^\infty \langle L_\theta(r)v(r,\cdot),v(r,\cdot)\rangle dr$$
$$-\int_{r_1}^\infty(E+\varphi'^2)\|v(r,\cdot)\|^2dr
+\int_{r_1}^\infty\langle (V+h^2q)v(r,\cdot),v(r,\cdot)\rangle dr.$$
Since $\varphi'={\cal O}(\tau)$, this implies
$$\int_{r_1}^\infty\|{\cal D}_rv(r,\cdot)\|^2dr\le C_1\tau^2\int_{r_1}^\infty\|v(r,\cdot)\|^2dr$$
$$+\left(\int_{r_1}^\infty r^{-2s}\|v(r,\cdot)\|^2dr\right)^{1/2}
\left(\int_{r_1}^\infty r^{2s}\|{\cal P}^\pm_\varphi(h)v(r,\cdot)\|^2dr\right)^{1/2}
$$
with some constant $C_1>0$. 
Hence
\begin{equation}\label{eq:4.7}
\varepsilon h^{-1}f(a)^2\int_{r_1}^\infty\|{\cal D}_rv(r,\cdot)\|^2dr\le C_1\tau^2\varepsilon h^{-1}f(a)^2\int_{r_1}^\infty\|v(r,\cdot)\|^2dr$$
$$+\gamma\epsilon\int_{r_1}^\infty r^{-2s}\|v(r,\cdot)\|^2dr
+\gamma^{-1}\epsilon^{-1}h^{-2}f(a)^4\int_{r_1}^\infty r^{2s}\|{\cal P}^\pm_\varphi(h)v(r,\cdot)\|^2dr
\end{equation}
for every $\gamma>0$. Taking $\gamma$ small enough, independent of $h$, and combining
the estimates (\ref{eq:4.6}) and (\ref{eq:4.7}), we get
\begin{equation}\label{eq:4.8}
\epsilon\int_{r_1}^\infty r^{-2s}\left(\|v(r,\cdot)\|^2+\|{\cal D}_rv(r,\cdot)\|^2\right)dr\le Cf(a)^4h^{-2}\epsilon^{-1}\int_{r_1}^\infty r^{2s}\|{\cal P}^\pm_\varphi(h)v(r,\cdot)\|^2dr$$ $$
+C\varepsilon h^{-1}f(a)^2\tau^2\int_{r_1}^\infty\|v(r,\cdot)\|^2dr
\end{equation}
with a new constant $C>0$ independent of $h$ and $\varepsilon$. It is an easy observation now that the estimate 
(\ref{eq:4.8}) implies (\ref{eq:4.1}).
\eproof

\section{Proof of Theorem 1.1}

In this section we will derive Theorem 1.1 from Theorems 2.1 and 4.1. Let $r_1$ be as above and fix 
$r_j$, $j=2,3,4$, such that $r_1<r_2<r_3<r_4$. Choose functions $\eta_1, \eta_2\in C^\infty(M)$ such that
$\eta_1=1$ in $M\setminus Y_{r_1}$, $\eta_1=0$ in $Y_{r_2}$, $\eta_2=1$ in $M\setminus Y_{r_3}$, $\eta_2=0$ in $Y_{r_4}$, 
and $\eta_j|_{Y_{r_1}}$ depending only on the variable $r$. Then we have
$$[P(h),\eta_j]=-h^2[\Delta_g,\eta_j]=-2h^2\eta'_j\partial_r-h^2\eta''_j-h^2\eta'_jp'p^{-1}.$$
Let $u\in H^2(M,d{\rm Vol}_g)$ be such that $\chi_s^{-1}(P(h)-E\pm i\varepsilon)u\in L^2(M,d{\rm Vol}_g)$. If $\partial M\neq\emptyset$
we require that $u|_{\partial M}=0$. 
Set 
$${\cal Q}_0=\|\chi_s^{-1}(P(h)-E\pm i\varepsilon)u\|_{L^2(M,d{\rm Vol}_g)},$$
$${\cal Q}_1=\|u\|_{L^2(Y_{r_1}\setminus Y_{r_2})}+\|{\cal D}_ru\|_{L^2(Y_{r_1}\setminus Y_{r_2})},$$
$${\cal Q}_2=\|u\|_{L^2(Y_{r_3}\setminus Y_{r_4})}+\|{\cal D}_ru\|_{L^2(Y_{r_3}\setminus Y_{r_4})},$$
and observe that
$$\|[P(h),\eta_j]u\|_{L^2}\lesssim {\cal Q}_j,\quad j=1,2.$$
We now apply Theorem 2.1 to the function $\eta_2u$ to obtain
\begin{equation}\label{eq:5.1}
\|u\|_{H_h^1(M\setminus Y_{r_3})}\le \|\eta_2u\|_{H_h^1(M\setminus Y_{r_4})}\le e^{\gamma h^{-4/3}}\|(P(h)-E\pm i\varepsilon)\eta_2u\|_{L^2(M\setminus Y_{r_4})}$$
 $$\le e^{\gamma h^{-4/3}}\|(P(h)-E\pm i\varepsilon)u\|_{L^2(M\setminus Y_{r_4})}+e^{\gamma h^{-4/3}}{\cal Q}_2
\end{equation}
with probably a new constant $\gamma>0$. 
In particular, (\ref{eq:5.1}) implies
\begin{equation}\label{eq:5.2}
{\cal Q}_1\le e^{\gamma h^{-4/3}}{\cal Q}_0+e^{\gamma h^{-4/3}}{\cal Q}_2.
\end{equation}
On the other hand, Theorem 4.1 applied to the function $(1-\eta_1)u$ yields 
\begin{equation}\label{eq:5.3}
\|r^{-s}e^{\varphi/h}u\|_{L^2(Y_{r_2},d{\rm Vol}_g)}+\|r^{-s}e^{\varphi/h}{\cal D}_ru\|_{L^2(Y_{r_2},d{\rm Vol}_g)}$$
 $$\le\|r^{-s}e^{\varphi/h}(1-\eta_1)u\|_{L^2(Y_{r_1},d{\rm Vol}_g)}+\|r^{-s}e^{\varphi/h}{\cal D}_r(1-\eta_1)u\|_{L^2(Y_{r_1},d{\rm Vol}_g)}$$
$$\le Cf(a)^2(\epsilon h)^{-1}\|r^se^{\varphi/h}(P(h)-E\pm i\varepsilon)(1-\eta_1)u\|_{L^2(Y_{r_1},d{\rm Vol}_g)}$$ 
$$+C\tau f(a)\varepsilon^{1/2}(\epsilon h)^{-1/2}\|e^{\varphi/h}u\|_{L^2(Y_{r_1},d{\rm Vol}_g)}$$
$$\le Cf(a)^2(\epsilon h)^{-1}\|r^se^{\varphi/h}(P(h)-E\pm i\varepsilon)u\|_{L^2(Y_{r_1},d{\rm Vol}_g)}+Cf(a)^2(\epsilon h)^{-1}e^{\varphi(r_2)/h}{\cal Q}_1$$ 
$$+C\tau f(a)\varepsilon^{1/2}(\epsilon h)^{-1/2}\|e^{\varphi/h}u\|_{L^2(Y_{r_1},d{\rm Vol}_g)}.
\end{equation}
In particular, (\ref{eq:5.3}) implies
\begin{equation}\label{eq:5.4}
e^{\varphi(r_3)/h}{\cal Q}_2\le Cf(a)^2(\epsilon h)^{-1}e^{\max\varphi/h}{\cal Q}_0+C\tau f(a)\varepsilon^{1/2}(\epsilon h)^{-1/2}
e^{\max\varphi/h}\|u\|_{L^2(M,d{\rm Vol}_g)}$$
$$+Cf(a)^2(\epsilon h)^{-1}e^{\varphi(r_2)/h}{\cal Q}_1.
\end{equation}
We have 
$$\varphi(r_3)-\varphi(r_2)=\tau\int_{r_2}^{r_3}\left(f(r)^{-1}-f(a)^{-1}\right)dr\ge c\tau_0h^{-1/3}$$
with some constant $c>0$. Observe also that $f(a)={\cal O}(h^{-km})$ if $f$ is given by (\ref{eq:1.3}), while in the other case the
assumption $\delta>\frac{3\alpha}{4}+1$ guarantees that $f(a)={\cal O}(e^{h^{-4/3}})$. Thus in both cases we deduce from (\ref{eq:5.4})
\begin{equation}\label{eq:5.5}
{\cal Q}_2\le \exp\left(\beta h^{-4/3}+\max\varphi/h\right){\cal Q}_0+\varepsilon^{1/2}
\exp\left(\beta h^{-4/3}+\max\varphi/h\right)\|u\|_{L^2(M,d{\rm Vol}_g)}$$
$$+\exp\left((\beta-c\tau_0)h^{-4/3}\right){\cal Q}_1
\end{equation}
with some constant $\beta>0$. Combining (\ref{eq:5.2}) and (\ref{eq:5.5}) we get
\begin{equation}\label{eq:5.6}
{\cal Q}_2\le \exp\left((\beta+\gamma)h^{-4/3}+\max\varphi/h\right){\cal Q}_0+\varepsilon^{1/2}
\exp\left(\beta h^{-4/3}+\max\varphi/h\right)\|u\|_{L^2(M,d{\rm Vol}_g)}$$
$$+\exp\left((\beta+\gamma-c\tau_0)h^{-4/3}\right){\cal Q}_2.
\end{equation}
Taking $\tau_0$ big enough and $h$ small enough, we can absorb the last term in the right-hand side of (\ref{eq:5.6}) 
to conclude that 
\begin{equation}\label{eq:5.7}
{\cal Q}_1+{\cal Q}_2\le \exp\left(\beta_1 h^{-4/3}+\max\varphi/h\right){\cal Q}_0+\varepsilon^{1/2}
\exp\left(\beta_1 h^{-4/3}+\max\varphi/h\right)\|u\|_{L^2(M,d{\rm Vol}_g)}
\end{equation}
with some constant $\beta_1>0$. By (\ref{eq:5.1}), (\ref{eq:5.3}) and (\ref{eq:5.7}) we obtain
\begin{equation}\label{eq:5.8}
\|\chi_su\|_{L^2(M,d{\rm Vol}_g)}\le N{\cal Q}_0+\varepsilon^{1/2}
N\|u\|_{L^2(M,d{\rm Vol}_g)}
\end{equation}
where
$$N=\exp\left(\beta_2 h^{-4/3}+\max\varphi/h\right)$$
with some constant $\beta_2>0$. On the other hand, since the operator $P(h)$ is symmetric, we have
\begin{equation}\label{eq:5.9}
\varepsilon\|u\|^2_{L^2(M,d{\rm Vol}_g)}=\pm{\rm Im}\,\langle (P(h)-E\pm i\varepsilon)u,u\rangle_{L^2(M,d{\rm Vol}_g)}$$
$$\le (2N)^{-2}\|\chi_su\|^2_{L^2(M,d{\rm Vol}_g)}+(2N)^2\|\chi_s^{-1}(P(h)-E\pm i\varepsilon)u\|^2_{L^2(M,d{\rm Vol}_g)}.
\end{equation}
We rewrite (\ref{eq:5.9}) in the form
\begin{equation}\label{eq:5.10}
N\varepsilon^{1/2}\|u\|_{L^2(M,d{\rm Vol}_g)}\le \frac{1}{2}\|\chi_su\|_{L^2(M,d{\rm Vol}_g)}+
2N^2\|\chi_s^{-1}(P(h)-E\pm i\varepsilon)u\|_{L^2(M,d{\rm Vol}_g)}.
\end{equation}
Combining (\ref{eq:5.8}) and (\ref{eq:5.10}) we get
\begin{equation}\label{eq:5.11}
\|\chi_su\|_{L^2(M,d{\rm Vol}_g)}\le 4N^2\|\chi_s^{-1}(P(h)-E\pm i\varepsilon)u\|_{L^2(M,d{\rm Vol}_g)}.
\end{equation}
It follows from (\ref{eq:5.11}) that the resolvent estimate
\begin{equation}\label{eq:5.12}
\left\|\chi_s(P(h)-E\pm i\varepsilon)^{-1}\chi_s
\right\|_{L^2(M,d{\rm Vol}_g)\to L^2(M,d{\rm Vol}_g)}\le 4N^2
\end{equation}
holds for all $0<h\ll 1$, $0<\varepsilon\le 1$ and $s$ satisfying (\ref{eq:3.1}).
 Observe also that if (\ref{eq:5.12}) holds for $s$ satisfying (\ref{eq:3.1}),
it holds for all $s>1/2$ independent of $h$. Thus, Theorem 1.1 follows from the bound (\ref{eq:5.12}) and Lemma 3.2.

\end{document}